\documentclass[a4paper,10pt]{article}
\usepackage[latin1]{inputenc}
\usepackage[T1]{fontenc}
\usepackage{amssymb}
\usepackage{epsfig}
\usepackage{amssymb}
\usepackage{amsmath}
\usepackage{amsfonts}
\usepackage{graphicx, url}
\usepackage{makeidx}
\usepackage{cite}
\usepackage{fancyhdr}
\usepackage{color}
\usepackage{setspace}

\newtheorem{theorem}{Theorem}[section]

\newtheorem{remark}[theorem]{Remark}

\newenvironment {proof} {{\it Proof.}}{\hspace*{\fill}$\Box$\par\vspace{4mm}}
\usepackage{amsmath}
\usepackage{amsfonts}

\newcommand{\mc}{\mathcal}

\title{A duopoly preemption game with two alternative stochastic investment choices}

\author{K. R. Dahl\footnotemark[2]\ \footnotemark[5]
\and E. Stokkereit\footnotemark[3]}

\begin{document}

\maketitle

\renewcommand{\thefootnote}{\fnsymbol{footnote}}

\footnotetext[2]{Department of Mathematics, University of Oslo, Pb. 1053 Blindern, 0316 Oslo, Norway. kristrd@math.uio.no}

\footnotetext[3]{Department of Economics, University of Oslo, Pb. 1095 Blindern, 0316 Oslo, Norway. espen.stokkereit@econ.uio.no}

\begin{abstract}
This paper studies a duopoly investment model with uncertainty. There are two alternative irreversible investments. The first firm to invest gets a monopoly benefit for a specified period of time. The second firm to invest gets information based on what happens with the first investor, as well as cost reduction benefits. We describe the payoff functions for both the leader and follower firm. Then, we present a stochastic control game where the firms can choose when to invest, and hence influence whether they become the leader or the follower. In order to solve this problem, we combine techniques from optimal stopping and game theory. For a specific choice of parametres, we show that no pure symmetric subgame perfect Nash equilibrium exists. However, an asymmetric equilibrium is characterized. In this equilibrium, two disjoint intervals of market demand level give rise to preemptive investment behavior of the firms, while the firms otherwise are more reluctant to be the first mover.
\end{abstract}

%\begin{keywords}

\medskip


\section{Introduction}

This paper studies how a two-firm preemptive investment game is affected when the firms involved have two investment opportunities instead of one. We compare this to the classic duopoly investment model, see e.g. Fudenberg and Tirole \cite{Fudenberg}. 

The analysis is motivated by D{\'e}camps et al~\cite{Decamps}. They show that adding another investment choice to a monopoly game gives rise to a dichotomous investment region. That is, the investment region is no longer a connected set, and it is this dichotomy that causes the preemption game to change structure.

We model a duopoly which  operates in a market where there is a stochastic, exogenous demand, driven by a geometric Brownian motion (GBM). The firms can at any time carry out a costly investment, where they must choose between two different technologies. The investment is irreversible, and can only be done one time by each firm. The investment reward is a higher profit rate and hence a greater income. However, how much the technologies increase the firm's profit is uncertain; each investment can either end up in a high profit state or a low profit state. As soon as the leader has invested, the profit state of his investment choice is revealed to both firms. Also, for the follower the investment cost of copying the same technology choice is reduced. Hence, being the second investor has benefits. On the other hand, the first mover will attain a monopoly benefit due to technological advantage until the second mover has invested. We solve this problem by using results by D{\'e}camps et al~\cite{Decamps} combined with game theoretic techniques based on the analysis of Chevalier-Roignant and Trigiorgis~\cite{TrigiorgisCompetitiveStrategies}.

\subsection{Literature overview}
Duopoly games often have complicated structures where the firms both compete and cooperate. An important example is the case of cooperative research effort, which is studied in Aspremont and Jacquemin~\cite{AspremontJacquemin}. The paper models this cooperation to be active, in the sense that the firms at a R\&D stage maximize their joint profit function, and then compete at a "market stage". Several other papers study this idea, but many exclude an active cooperative setting, and instead model the benefits of other firms' investment by technology spillovers. Technology spillover is the concept that when one firm attains a new technology, it becomes easier for other firms to acquire the same technology, for instance because they can learn from the innovator firm's research advances. An example of such a paper is Femminis and Martini~\cite{FemminisMartini}. In this paper, investment is a "one-shot" action which increases the firms' profit rate. A certain time after the first firm has invested, the investment cost is reduced for the second firm. Thus, the first firm can be seen as an innovator, while the follower can be seen as an imitator.% By some results related to GBM, the problem reduces from a supremum over integral to a more simple expression.
%However, one may argue that a weakness of Femminis and Martini's~\cite{FemminisMartini} paper is that the stochastic component is not connected to the pay-off of investment, but instead to the general profit rate. Also, since the choice is restricted to invest or not invest, and this investment is one-shot, the model is not suitable for deciding long run behavior. For instance, there is never a possibility for the second mover to change into a leading innovator.
%This is similar to our setting, however, in comparison to Femminis and Martini, we introduce an extra source of uncertainty by making the increase of the profit rate random. This is also in contrast to the uncertainty modeled by Dasgupta and Stiglitz~\cite{DasguptaStiglitz}.%, where investment is uncertain because after investing the time (and not size) of payoff is uncertain.

%First, what if investment was not simply a one time shot, but that the first investment leads to the possibility of a second investment etc. An "investment chain". This could be finite or countable infinite. I still think the problem would be reduced from an integral as before, but a little bit more complicated.
%Here, there are several ways to do things. For instance, the first mover could make it cheaper to invest for the second mover. Or, if for instance the first mover are $n$ investment steps ahead, the second mover can skip a certain amount of investment of investment steps ($\leq n$).

Some papers that focus on second mover advantages, due to for instance information benefits, are Hoppe~\cite{Hoppe} and Y. J. Yap, S. Luckraz and S. K. Tey~\cite{YapEtAl}. A paper with an in-depth study of imitation versus innovation is Bessen and Maskin~\cite{BessenMaskin}. They consider a setting where the more firms that invest, the bigger the probability of innovating becomes. However, if one firm succeeds, the other firms can copy this. This gives incentives to not invest and instead simply imitate. The paper focuses on social optimality, and compares situations with and without patents. Also, they cover both a static and a sequential investment context. Another relevant paper that examines the choice of doing innovation versus imitation is Jin and Troege~\cite{JinTroege}.

Huisman and Kort~\cite{HuismanKort} analyzes a choice between two investments. However, the setting of this paper is different from ours, as one does not know when the second investment opportunity arrives. Other papers studying investment choice are D{\'e}camps et al.~\cite{Decamps} and Nishihara and Ohyama~\cite{Nishihara}. The solution method we apply in this paper is based on the results in D{\'e}camps et al.~\cite{Decamps}. However, in contrast to Decamps et al.~\cite{Decamps}, we study a game between two firms and we have an additional uncertainty in the profit rate obtained from the investment. Nishihara and Ohyama~\cite{Nishihara} studies a game, but without uncertain profit rate.

The papers by Huisman et al.~\cite{Thijssen_information} and \cite{Thijssen} both analyze stochastic games related to real options models. The paper~\cite{Thijssen} is quite general, but the firms only have one investment choice. We extend this by adding another technology. They study a preemption example, and we extend this example in a way which results in a more complicated structure of behavior of the firms (contrast Figure 1 in Huisman et al.~\cite{Thijssen} to our Figure \ref{fig:awesome_image}). The paper Huisman et al.~\cite{Thijssen_information} differs from ours in that they don't have a stochastic demand, but rather a Poisson process for when new information is revealed to the firms.

%In the mathematical finance literature, a game where it is predetermined which player goes first is called a Stackelberg game. This has, among others, been studied by Bensoussan et al. in ~\cite{BensoussanEtAl} and \cite{BensoussanEtAl2} and by {\O}ksendal et al. in ~\cite{OksendalUboe} and \cite{OksendalUboe2}. Our initial problem is of this form, however we relax this assumption and aim to create a dynamic game in the final sections of the paper by introducing a patent. The setting of the Stackelberg game presented by Bensoussan et al.~\cite{BensoussanEtAl2} is similar to our initial setting. However, \cite{BensoussanEtAl2} consider a demand process which is an arithmetic Brownian motion, while we consider a geometric Brownian motion. Also, we use different methods to solve the problem than \cite{BensoussanEtAl2}.

\subsection{Structure of the paper}

The paper is organized as follows. In Section \ref{The model}, we give a mathematical formulation of the model. Then, Section \ref {sec: symmetric} describes the different payoff functions the firms may attain (depending on the investment outcome). This includes the case of a simultaneous Cournot investment, being the market follower and being the market leader. Section \ref{Example_preemption} presents a stochastic game in this framework. The game is defined in Section \ref{definingthegame}, based on the approach of Huisman~\cite{Thijssen}. Then a specific scenario where double preemption intervals occur is analyzed in Section \ref{specialscenario}. We show that no subgame perfect symmetric Nash equilibrium exists in this case, however, we characterize an asymmetric equilibrium. Finally, we conclude in Section \ref{conclusion}.

\section{The model}
\label{The model}

Consider a framework with two firms, both facing two different irreversible investment opportunities. We call these investment 1 and investment 2. Both firms may invest at any time in either of the investment opportunities, but not in both. When the first firm has invested, thus becoming the leader, nature chooses one of two outcomes: a high profit outcome or a low profit outcome. This outcome is determined according to some probability measure. Hence, we consider a finite probability space $(\Omega_F, \mc{F}_F, P_F)$ (the subscript $F$ stands for ``finite''), where $\Omega_F = \{\omega_l^l, \omega_l^h, \omega_h^l, \omega_h^h\}$, $\mc{F}_F$ is the corresponding $\sigma$-algebra and $P_F$ is a probability measure. Here, $\omega_h^l$ is the scenario where investing in alternative 1 leads to a high profit, denoted $\pi_{1}^h>0$, while investment 2 leads to low profit, $\pi_{2}^l>0$. Similar interpretations hold for $\omega_l^l, \omega_l^h$ and $\omega_h^h$. These probabilities are the same for both firms.

\begin{remark}
The motivation for this problem is to study a model for an investment game involving uncertainty, since this is more realistic than a deterministic model. Only considering two possible outcomes is a simplification which is made to make the computations more manageable. However, similar kinds of arguments can be made for some number $n < \infty$ of potential outcomes, though this will become chaotic as $n$ grows.
\end{remark}

The second firm to invest, thus becoming the follower, can choose what to do after the leader has invested. The follower observes the outcome of the leader's investment. Hence, the follower has the advantage of having more information than the leader.

The profit rates of the two firms vary according to the state of the game. In accordance to what is often done in the literature, we normalize such that both agents get no profit before any investment has occurred. When the first firm has invested, it immediately gets either profit rate $\pi_{1}^{h}$ (high profit, investment 1), $\pi_{1}^{l}$ (low profit, investment $1$), $\pi_{2}^{h}$ or $\pi_{2}^{l}$, depending on whether the firm chooses investment 1 or 2, and on which scenario $\omega \in \Omega_F$ is realized. Until the second firm invests, the leader also enjoys a monopoly benefit $\xi$. The follower continues to have profit rate $0$ until it invests. If it chooses to copy the technology of the first investor, by for example choosing investment 1, the follower gets the same level of profit as the first investor (either $\pi_{1}^{h}$ or $\pi_{1}^{l}$ depending on the outcome of the $\Omega_F$-draw). However, if the second investor chooses not to follow the first investor (i.e. by choosing investment 2), it will get a random profit rate, either $\pi_{2}^{h}$ or $\pi_{2}^{l}$, according to the probability measure $P_F$ given that only $\pi_{2}^{h}$ or $\pi_{2}^{l}$ can happen.

To simplify, we assume that the problem is symmetric, i.e. that $\pi_1^h=\pi_2^h =: \pi_h$, $\pi_1^l=\pi_2^l =: \pi_l$ and that the probabilities of high or low outcomes, respectively, are the same for both investments. This means that we can remove the technology choice for the first firm, and just look at nature's choice of high or low profit.

%ILLUSTRASJON AV SCENARIOTRE:
%
 \begin{figure}[h]
   \setlength{\unitlength}{0.7mm}
 \begin{picture}(100,100)(-40,0)

  %Dots
  \put(20,45){\circle*{3}} % t=1

  \put(40,65){\circle*{3}} % t=2
  \put(40,25){\circle*{3}} % t=2

  \put(60,13){\circle*{3}} % t=3
  \put(60,25){\circle*{3}} % t=3
  \put(60,65){\circle*{3}} % t=3
  \put(60,76){\circle*{3}} % t=3

  \put(80,65){\circle*{3}} % t=4
  \put(80,76){\circle*{3}} % t=4
  \put(80,53){\circle*{3}} % t=4
  %\put(80,88){\circle*{3}} % t=4
  \put(80,13){\circle*{3}} % t=4
  \put(80,25){\circle*{3}} % t=4
  \put(80,2){\circle*{3}} % t=4
  %\put(80,36){\circle*{3}} % t=4

  %Lines
  \put(20,45){\line(1,1){20}}
  \put(20,45){\line(1,-1){20}}

  \put(40,65){\line(5,3){20}}
  \put(40,65){\line(5,0){20}}
  \put(40,25){\line(5,0){20}}
  \put(40,25){\line(5,-3){20}}

  \put(60,65){\line(5,0){20}}
  \put(60,65){\line(5,-3){20}}
  \put(60,76){\line(5,0){20}}
  %\put(60,76){\line(5,3){20}}
  \put(60,13){\line(5,0){20}}
  \put(60,13){\line(5,-3){20}}
  \put(60,25){\line(5,0){20}}
%  \put(60,25){\line(5,3){20}}

%Legger til tekst pa nodene
  \put(2,53){\makebox(0,0){1. investor}}
  \put(36,70){\makebox(0,0){$\pi_h$}}
  \put(36,20){\makebox(0,0){$\pi_l$}}
  \put(34,46){\makebox(0,0){Nature}}
  %\put(58,82){\makebox(0,0){$\pi_1^h$}}
  %\put(58,60){\makebox(0,0){$\pi_1^l$}}
  %\put(58,30){\makebox(0,0){$\pi_2^h$}}
  %\put(58,8){\makebox(0,0){$\pi_2^l$}}
  \put(47,95){\makebox(0,0){2. investor}}

  \put(48,75){\makebox(0,0){\tiny{Follow}}}
  \put(51,60){\makebox(0,0){\tiny{Innovate}}}
  \put(51,28){\makebox(0,0){\tiny{Follow}}}
  \put(48,15){\makebox(0,0){\tiny{Innovate}}}

  \put(75,44){\makebox(0,0){Nature}}
  \put(88,1){\makebox(0,0){$\pi_l$}}
  \put(88,12){\makebox(0,0){$\pi_h$}}
  \put(88,25){\makebox(0,0){$\pi_l$}}
  \put(88,53){\makebox(0,0){$\pi_l$}}
  \put(88,64){\makebox(0,0){$\pi_h$}}
  \put(88,75){\makebox(0,0){$\pi_h$}}

 \end{picture}
 \caption{A decision tree.}
 \label{fig: decisiontree}
 \end{figure}
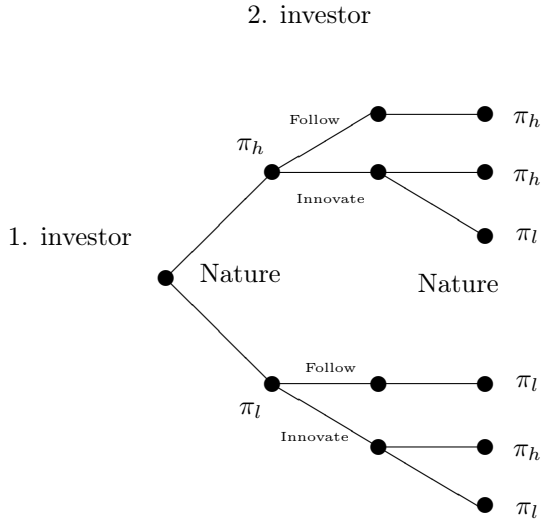

To model the exogenous and uncertain demand for the duopoly's product, we introduce another probability space $(\Omega, \mc{F}, P)$, and consider a Brownian motion $(B_t)_{t \in [0, \infty)}$ in this space. Let $(\mc{F}_t)_{t \in [0,\infty)}$ be the filtration generated by the Brownian motion. The process modeling the uncertainty in demand is a geometric Brownian motion which we denote by $(Z_t)_{t \in [0, \infty)}$, such that

\[
 dZ_t = \alpha Z_t dt + \sigma Z_t dB_t.
\]

\noindent Here, $\alpha, \sigma \in \mathbb{R}_{+}$. Also, let $r \in \mathbb{R}_+$ be the interest rate. We assume $\alpha < r$.

%\textcolor{magenta}{(STEMMER DETTE? VIRKER LITT RART AT DET SIKRE MAN FAAR VIA RENTA SKAL VAERE MER ENN DET MAN FAAR VED AA GJOERE EN USIKKER INVESTERING, SELV OM MAN JO KAN TJENE PAA DEN USIKRE INVESTERINGEN. DET ER DETTE FEMMINIS OG MARTINI ANTAR, MEN SYNES DU DET VIRKER RIMELIG?)}.

For the first firm to invest, we let the cost of any investment be $\mc{I} > 0$. If the second investor chooses to follow the leader, the investment cost is reduced to $(1-\theta)\mc{I}$, $\theta \in (0,1)$, while the price for the other investment choice is still $\mc{I}$. The intuition is that it is more expensive to be innovative than to imitate technology that already exists. After both firms have invested, no more choices are made and the profit rates of the investors stay constant.

\section{The payoff functions}
\label{sec: symmetric}

%We assume that the problem is symmetric in the sense that $\pi_{1}^h = \pi_{2}^h$, $\pi_{1}^l = \pi_{2}^l$ and
%
%\[
%P_F(\cdot | \mbox{ player I chooses investment 1}) = P_F(\cdot | \mbox{ player I chooses investment 2}).
%\]

%Let $\pi_{1, j}$, $j \in \{A,B\}$ be the $\Omega_F$-random variable which takes value $\pi_{1,j}(\omega) = \pi_{1,j}^h$ for $\omega \in \{\omega_h^h, \omega_h^l\}$ and $\pi_{1,j}(\omega) = \pi_{1,j}^l$ for $\omega \in \{\omega_l^h, \omega_l^l\}$, $j \in \{A,B\}$. Also, let the $\Omega_F$-random variable $\pi_{2,j}$ be defined similarly.

%Furthermore, we define $p := E[\bb{1}_A]= P(A)$. Thus, $p$ is the probability that the follower gets high profit when choosing another investment path than the leader.

We want to describe a game where the firms can choose when to invest, so they can influence whether they become the leader or the follower. In order to do this, we must compute the expected payoff functions attained for the different roles. More specifically, we will compute the following functions for all times $t \geq 0$, where we assume that the stochastic demand process $Z_t$ has value $z_t$ at time $t$:

\begin{enumerate}
 \item[$(i)$] $C(z_t)$: The payoff of both firms when they invest simultaneously at time $t$, and the firms do not take the roles of leader and follower. Note that the $C$ stands for ``Cournot'' (see also Chevalier-Roignant and Trigeorgis~\cite{TrigiorgisCompetitiveStrategies}). In this case, we assume that neither of the firms get any information advantage or reduced investment cost.
 \item[$(ii)$] $F(z_t)$: The payoff of the second investor given that the leader stops at time t and that the follower is acting optimally.
 \item[$(iii)$] $L(z_t)$: The payoff of the first investor who stops at time t and hence becomes the leader, given that the other firm will act optimally as the follower.
\end{enumerate}

We divide the computation of these functions into several subsections.

In the following, $E[\cdot]$ denotes the expectation with respect to the product measure $P \times P_F$, also let $E^P[\cdot]$ and $E^{P_F}[\cdot]$ denote the expectation w.r.t. the measures $P$ and $P_F$ respectively.

\subsection{The Cournot case: $C(z_t)$}

Let $\pi$ be the random variable taking the values $\pi_l$ and $\pi_h$ with probabilities determined by $P_F$. Let $E^{\cdot}_{z_t}[\cdot]$ denote the expectation given that $Z_t=z_t$.

In the case where both firms invest exactly at time $t$ and the Cournot outcome is attained, their payoff functions are identical, and given by
\begin{equation}
 \label{eq: Cournot}
 \begin{array}{lll}
  C(z_t) &=& E_{z_t}[\int_t^{\infty} Z(s) e^{-r(s-t)} \pi ds - \mc{I}] \\[\smallskipamount]
      &=&  E_{z_t}[\int_0^{\infty} Z(u+t) e^{-ru} E[\pi] du - \mc{I}] \\[\smallskipamount]
      &=& E^{P_F}[\pi] \int_0^{\infty} E^P_{z_t}[Z(u+t)] e^{-ru} du - \mc{I} \\[\smallskipamount]
      &=& E^{P_F}[\pi]\frac{z_t}{r-\alpha} - \mc{I}.
 \end{array}
\end{equation}
\noindent Here, the first and the third equalities follow from Tonelli's theorem, the second equality from a change of variables and the final equality from that $Z$ is a geometric Brownian motion.

\subsection{The follower's payoff: $F(z_t)$}

We begin by finding the payoff function of the follower, hence we assume that one of the firms has already invested at time $t$. Based on this we can later determine the payoff of the leader who's optimal behavior depends on the behavior of the follower.

For the follower, two scenarios can occur: The leader got high profit $\pi_h$ or low profit $\pi_l$. The follower knows which profit level the leader got, and can choose to copy or innovate based on this.

\subsubsection{The case where the leader got low profit: $F_L(z_t)$}
\label{sec: F_L}

If the leader got the low profit $\pi_l$, the payoff function of the follower at the leader's investment time $t$ is:
\begin{equation}
\label{eq: follower_low}
\begin{array}{lll}
  F_L(z_t) := e^{r t} \sup_{\tau \geq 0 } \max \Big\{ E^P \big[ \{ \int_{\tau+t}^{\infty} \pi_{l} Z_s e^{-r s} ds -(1-\theta)\mc{I} e^{-r (\tau+ t)} \big], \\[\smallskipamount]
   \hspace{2cm}  E^P \big[ \int_{\tau+t}^{\infty}  E^{P_F}[\pi] Z_s e^{-r s} ds  - \mc{I} e^{-r (\tau + t)} \big] \Big\}
 \end{array}
\end{equation}
\noindent where $\tau$ is the stopping time from the investment time of the leader until the investment time of the follower (that is, $\tau + t$ is the investment time of the follower).

The first expression in the maximum corresponds to copying the leader, while the second expression corresponds to innovating. The first expression can be rewritten:
\[
 \begin{array}{lll}
  E^P \Big[ \{ \int_{\tau + t}^{\infty} \pi_l Z_s e^{-r s} ds -(1-\theta)\mc{I} e^{-r (\tau + t)} \Big] \\[\smallskipamount]
  \hspace{3cm} = e^{-rt} E^P _{z_t} \Big[ \int_{\tau}^{\infty} \pi_l Z_s e^{-rs} ds -(1-\theta)\mc{I} e^{-r \tau } \Big] \\[\smallskipamount]
  \hspace{3cm} = e^{-rt} E^P_{z_t} \Big[ E^P \big[ \int_{\tau }^{\infty} \pi_{l} Z_s e^{-r s} ds -(1-\theta)\mc{I} e^{-r \tau } | \mc{F}_{\tau} \Big] \\[\medskipamount]
  %
  %\hspace{3cm} = e^{-rt} E_{z_t} \Big[ e^{-r \tau} E_{Z(\tau)} \big[ \int_{0}^{\infty} \pi_l Z_s e^{-r s} ds -(1-\theta)\mc{I} \big]  \Big] \\[\medskipamount]
  %
  \hspace{3cm} = e^{-rt} E^P_{z_t} \Big[ e^{-r \tau} \big( \int_{0}^{\infty} \pi_{l} Z(\tau) e^{(\alpha-r) s} ds -(1-\theta)\mc{I} \big) \Big] \\[\medskipamount]
  \hspace{3cm} =e^{-rt} E^P_{z_t} \Big[ e^{-r \tau} \frac{Z(\tau) \pi_{l}}{r-\alpha} -(1-\theta)\mc{I}e^{-r\tau} \Big]
 \end{array}
\]

\noindent where the second equality follows from double expectation, the third equality from the strong Markov property for the It{\^o} diffusion $Z_t$ (see also {\O}ksendal~\cite{Oksendal} exercise 7.11 for justification of the integral term) and the final equality follows from that $Z_t$ is a geometric Brownian motion.

A similar computation shows that the second part of the maximum in problem~\eqref{eq: follower_low} can be written

\[
\begin{array}{lll}
  E_{z_t}^{P} \big[ \int_{\tau +t}^{\infty}  E^{P_F}[\pi] Z_s e^{-r s} ds  - \mc{I} e^{-r (\tau + t)} \big] \\[\smallskipamount]
  \hspace{2cm} =e^{-rt} E_{z_t}^{P} \Big[ e^{-r \tau} \frac{Z(\tau) E^{P_F}[\pi]}{r-\alpha} -\mc{I}e^{-r\tau} \Big].

\end{array}
\]

Hence, problem~\eqref{eq: follower_low} is equivalent to

\begin{equation}
 \label{eq: follower_low2}
 F_L(z_t) = \sup_{\tau \geq 0} E_{z_t}^{P}[e^{-r \tau} \max_{i=1,2} \{ a_i Z(\tau) - K_i\}]
\end{equation}
\noindent where $a_1 := \frac{\pi_l}{r-\alpha}$, $a_2 := \frac{E[\pi]}{r-\alpha}$, $K_1 := (1-\theta) \mc{I} $ and $K_2 := \mc{I} $. Note that $K_1 < K_2$ and $a_1 < a_2$.

Problem~\eqref{eq: follower_low2} is of the same form as the problems in Decamps et al~\cite{Decamps} and Nishihara and Ohyama~\cite{Nishihara}. As in these papers, we have to consider several cases. The first case is $0< \frac{a_2}{a_1} <1$, however since $a_1 < a_2$, this never occurs.

Let $\gamma := \frac{1}{2} - \frac{\alpha}{\sigma^2} + \sqrt{( \frac{\alpha}{\sigma^2} - \frac{1}{2})^2 + \frac{2r}{\sigma^2}} > 1$. The next case to consider is $1< (\frac{a_2}{a_1})^{\frac{\gamma}{\gamma - 1}} <\frac{K_2}{K_1}$.  Then,

\begin{equation}
 \label{eq: follower_low_second_case}
 F_L(z_t) = \left\{
	\begin{array}{ll}
		A_0 z_t^{\gamma}  & \mbox{if } 0<z_t<z_{1,L}^* \\
		a_1 z_t - K_1 & \mbox{if } z_{1,L}^* \leq z_t \leq z_{2,L}^* \\
		B_0 z_t^{\gamma} - C_0 z_t^{\beta} & \mbox{if } z_{2,L}^* < z_t < z_{3,L}^* \\
		a_2 z_t - K_2 & \mbox{if } z_t \geq z_{3,L}^* \\
	\end{array}
\right.
\end{equation}
\noindent where $\beta := \frac{1}{2} - \frac{\alpha}{\sigma^2} - \sqrt{( \frac{\alpha}{\sigma^2} - \frac{1}{2})^2 + \frac{2r}{\sigma^2}} < 0$ and $A_0, B_0, C_0$ as well as the investment thresholds $z_{1,L}^*, z_{2,L}^*, z_{3,L}^*$ are determined using value matching and smooth pasting conditions as in Dixit and Pindyck~\cite{Dixit} (i.e. requiring continuity and continuous differentiability of the value function $F_L$). In particular, we find that

\begin{equation}
 \label{eq: z_{1,L}^*}
\begin{array}{lll}
 z_{1,L}^* &=& \frac{\gamma}{\gamma -1} \frac{(1-\theta) \mc{I} (r-\alpha)}{\pi_l}, \\[\smallskipamount]
 A_0 &=& \frac{(1-\theta) \mc{I}}{(\gamma -1)(z_{1,L}^*)^{\gamma}}.
 \end{array}
\end{equation}

Similarly, we can also find $z_{2,L}^*, z_{3,L}^*, B_0$ and $C_0$, however we omit writing these out since the expressions are long, and not necessary to know for our purposes. However, we will use the fact that $z_{1,L}^* \leq z_{2,L}^* \leq z_{3,L}^*$.

The values of this split function correspond to the following actions of the follower (starting at the top): 1) Wait until the demand reaches $z_{1,L}^*$, then copy the leader. 2) Copy right away. 3) Wait, but make no decision to copy or innovate. This decision depends on whether demand hits $z_{2,L}^*$ (copy) or $z_{3,L}^*$ (innovate) first. 4) Innovate right away.

The optimal stopping time is

\[
\tau^* := \inf \{ t \geq 0 | Z(t) \in [z_{1,L}^*, z_{2,L}^*] \cup [z_{3,L}^*, \infty) \}.
\]

The reason for the waiting behavior in the interval $(z_{2,L}^*, z_{3,L}^*)$ is due to the structure of the follower's payoff function
\[
 F_L(z_t) = \sup_{\tau \geq 0} E_{z_t}^{P}[e^{-r \tau} \max_{i=1,2} \{ a_i Z(\tau) - K_i\}],
\]
as shown in equation (\ref{eq: follower_low2}). Since $a_1 < a_2$ and $K_1 < K_2$ there is a point $\hat{z}$ where $a_1 \hat{z} - K_1 = a_2 \hat{z} - K_2 $, and at this point we have
\[
D^{+}F_L(\hat{z}) > D^{-}F_L(\hat{z}),
\]
where $D^{+}, D^{-}$ denotes the right and left derivatives respectively. Thus, at $\hat{z}$ the profit gain of a marginal increase in market demand is larger than the profit loss of a marginal decrease in market demand. As a result, there is an open nonempty interval containing $\hat{z}$ where the follower is interested in waiting to see how the market demand $Z$ evolves. For more on this issue, see the paper by Decamps et al~\cite{Decamps}.
\smallskip

The final case to consider is $(\frac{a_2}{a_1})^{\frac{\gamma}{\gamma - 1}} \geq \frac{K_2}{K_1}$. In this case,

\begin{equation}
 \label{eq: follower_low_third_case}
 F_L(z_t) = \left\{
	\begin{array}{ll}
		B_0 z_t^{\gamma}  & \mbox{if } 0<z_t<z_{3,L}^* \\
		a_2(t)z_t - K_2(t) & \mbox{if } z_t \geq z_{3,L}^*.
	\end{array}
\right.
\end{equation}

The optimal stopping time is $\tau^* = \inf \{ t \geq 0 | Z(t) \geq z_{3,L}^*\}$. In this situation, the follower waits until the demand reaches the level $z_{3,L}^*$ and then innovates.

\subsubsection{The case where the leader got high profit: $F_H (z_t)$}
\label{sec: F_H}

In the case that the leader got the high profit $\pi_h$, the payoff function of the follower at time $t$ has a similar form as the above situation:

\begin{equation}
\label{eq: follower_high}
\begin{array}{lll}
  F_H(z_t) := e^{r t} \sup_{\tau \geq 0} \max \Big\{ E^P \big[ \int_{\tau + t}^{\infty} \pi_{h} Z_s e^{-r s} ds -(1-\theta)\mc{I} e^{-r (\tau + t)} \big], \\[\smallskipamount]
   \hspace{2cm}  E^P \big[ \int_{\tau +t}^{\infty}  E^{P_F}[\pi] Z_s e^{-r s} ds  - \mc{I} e^{-r (\tau + t)} \big] \Big\}.
 \end{array}
\end{equation}
\noindent However, since both $\pi_h > E[\pi]$ and $(1 - \theta)\mc{I} < \mc{I}$ are true, it follows that the first term in the maximum expression always is higher than the second. Thus, the payoff function reduces to

\begin{equation}
 \label{eq: follower_high_2}
 F_H(z_t) = \left\{
	\begin{array}{ll}
		\frac{\mc{I}(1 - \theta)}{\gamma - 1} \frac{z_t}{z^*_H}^{\gamma}  & \mbox{if } 0<z_t< z^*_H \\[\smallskipamount]
		\frac{\pi_{h}}{r - \alpha}z_{t} - \mc{I}(1 - \theta) & \mbox{if } z_t \geq z^*_H.
	\end{array}
\right.
\end{equation}

Here,
\[
z^*_H = \frac{\gamma}{\gamma - 1}\frac{r - \alpha}{\pi_{h}}(1 - \theta)\mc{I}
\]
 is the critical value such that the follower waits for all $z_t < z^*_H$ and copies immediately for $z_t \geq z^*_H$ (see e.g. Femminis and Martini~\cite{FemminisMartini} for solution procedure), where $\gamma$ is as before.

 Note that $z^*_H < z^*_{1,L}$, i.e. the demand-threshold for the follower-firm investing when they are guaranteed a high profit is lower than the first investment threshold in the low profit case.

\subsection{The leader's payoff: $L(z_t)$}

Now that we have determined the optimal behavior and payoff function for the follower, we can derive the leader's optimal value function $L(z_t)$. This payoff function is the expectation (with respect to the probability measure corresponding to nature's choice of high or low profit) of the payoff functions in the case where the leader gets high or low profit respectively. We let $L_{L}(z_t)$ and $L_{H}(z_t)$ be the value functions in the low and high profit cases respectively. Then,

\begin{equation}
 \label{eq: leader_payoff}
 L(z_t)=P_F(\mbox{low profit}) L_{L}(z_t) + P_F(\mbox{high profit}) L_{H}(z_t).
\end{equation}

Since we have assumed that the rest of the problem is symmetric, we assume that
\[
P_F(\mbox{high profit}) = P_F(\mbox{low profit})=\frac{1}{2}.
\]

We derive the two value functions $L_{L}(z_t)$ and $L_{H}(z_t)$ in the following subsections.

\subsubsection{The low profit case: $L_{L}(z_t)$}
\label{lowprofcase}

The leader's payoff function depends on the behavior of the follower. Hence, we have to take the structure of the function $F_L(z_t)$ into consideration, see Section~\ref{sec: F_L}. As in Section~\ref{sec: F_L}, we consider two cases.

First, we consider the case where  $1< (\frac{a_2}{a_1})^{\frac{\gamma}{\gamma - 1}} <\frac{K_2}{K_1}$ (where the functions $a_1, a_2, K_1, K_2$ are defined as in Section~\ref{sec: F_L}). The payoff function $L_L(z_t)$ is the payoff of the leading firm when it invests at time $t$, given that the follower then behaves optimally.

By computations similar to those of Section~\ref{sec: F_L}, we can compute $L_L(z_t)$ (as we have assumed, due to symmetry, that the leader has no investment choice, and the investment time is $t$ per definition of $L_L(z_t)$). For example, in the case where $0<z_t<z_{1,L}^*$, we have
\[
 \begin{array}{lll}
  L_L(z_t) &=& E^P[\int_t^{t+\tau^*} (\pi_l + \xi) Z(s) e^{-r(s-t)} ds + \int_{\tau^* + t}^{\infty} \pi_l Z(s) e^{-r(s-t)}ds - \mc{I}] \\[\smallskipamount]

   &=&  E^P[\int_t^{\infty} (\pi_l + \xi) Z(s) e^{-r(s-t)} ds - \int_{\tau^*+t}^{\infty} \xi Z(s) e^{-r(s-t)}ds - \mc{I}] \\[\smallskipamount]

   &=& \frac{z_t( \pi_l + \xi) }{r-\alpha} - \frac{z_{1,L}^* \xi}{r-\alpha} E^P[e^{-r \tau^*}| Z(t)=z_t] - \mc{I} \\[\smallskipamount]

   &=& \frac{z_t( \pi_l + \xi) }{r-\alpha} - \frac{z_{1,L}^* \xi }{r-\alpha} e^{\nu \alpha_1(t) - |\alpha_1(t)| \sqrt{\nu^2 + 2r} } - \mc{I}\\[\smallskipamount]

   &=& \frac{\pi_l}{r-\alpha} z_t + \frac{\xi}{r-\alpha}z_t (1 - (\frac{z_t}{z_{1,L}^*})^{\gamma - 1} ) - \mc{I}
 \end{array}
\]

\noindent where $\alpha_1(t) := \frac{1}{\sigma} \ln(\frac{z_{1,L}^*}{z_t})$, $\nu := \frac{\alpha}{\sigma} - \frac{\sigma}{2}$. Note that in the first equality, the first integral corresponds to the monopoly benefit, which only lasts until the follower invests (at time $t+\tau^*$). Note also that the second equality follows by adding and subtracting the monopoly benefit, the third equality follows from a similar argument (using the strong Markov property for It{\^o} diffusions) as in Section~\ref{sec: F_L} and the fourth equality follows from that $\{Z_t\}_{t \geq 0}$ is a geometric Brownian motion, so we know the expectation of the first hitting time at the investment level (for the follower) $z_{1,L}^*$ (see also Jeanblanc~\cite{Jeanblanc}, section 8.1.3). The final equality follows from some basic algebra.

Thus, we have
\begin{equation}
 \label{eq: leader_low_second_case}
 L_L(z_t) = \left\{
	\begin{array}{ll}
		 \frac{\pi_l}{r-\alpha} z_t + \frac{\xi}{r-\alpha}z_t (1 - (\frac{z_t}{z_{1,L}^*})^{\gamma - 1} ) - \mc{I}  & \mbox{if } 0<z_t<z_{1,L}^* \\[\medskipamount]
		
		\frac{\pi_l }{r - \alpha}z_t - \mc{I} & \mbox{if } z_{1,L}^* \leq z_t \leq z_{2,L}^* \\[\medskipamount]
		
		\frac{\pi_l}{r-\alpha} z_t + M(z_t) - \mc{I} & \mbox{if } z_{2,L}^* \leq z_t \leq z_{3,L}^* \\[\medskipamount]
		
		 \frac{ \pi_l }{r - \alpha}z_t - \mc{I} & \mbox{if } z_t \geq z_{3,L}^* \\[\medskipamount]
	\end{array}
\right.
\end{equation}
%\noindent where $\alpha_i(t) := \frac{1}{\sigma} \ln(\frac{z_i^*}{z_t})$, $i=1,3$ (recall that $\nu = \frac{\alpha}{\sigma} - \frac{\sigma}{2}$), and

\noindent where,

\begin{equation}
\begin{array}{lll}
M(z_t) := (1- P(Z \mbox{ hits } z_{3,L}^* \mbox{ first} | Z(t)=z_t)) E^P[\int_t^{t+\tau^*}\xi Z(s) e^{-r(s-t)} ds| Z \mbox{ hits } z_{2,L}^* \mbox{ first}] \\

\hspace{0.4cm}+ P(Z \mbox{ hits } z_{3,L}^* \mbox{ first} | Z(t)=z_t) E^P[\int_t^{t+\tau^*}\xi Z(s) e^{-r(s-t)} ds| Z \mbox{ hits } z_{3,L}^* \mbox{ first}] .
\end{array}
\end{equation}

%\begin{equation}
%\begin{array}{lll}
%M(z_t) :=(1- P(Z \mbox{ hits } z_{3,L}^* \mbox{ first} | Z(t)=z_t)) \Big[ \frac{\pi_l}{r-\alpha} z_t + \frac{\xi}{r-\alpha}z_t (1 - (\frac{z_t}{z_{2,L}^*})^{\gamma - 1} ) \Big] \\
%
%\hspace{0.4cm}+ P(Z \mbox{ hits } z_{3,L}^* \mbox{ first} | Z(t)=z_t) \Big[ \frac{\pi_l}{r-\alpha} z_t + \frac{\xi}{r-\alpha}z_t (1 - (\frac{z_t}{z_{3,L}^*})^{\gamma - 1} ) \Big] .
%\end{array}
%\end{equation}

Here, using that $Z$ is a geometric Brownian motion (see Sigman~\cite{Sigman})

\[
 P(Z \mbox{ hits } z_{3,L}^* \mbox{ first} | Z(t)=z_t) = \frac{1 - e^{\tilde{\alpha}b}}{ e^{\tilde{\alpha}a} - e^{\tilde{\alpha}b}}
\]
where $a := \ln(\frac{z_{3,L}^*}{z_t})$, $b := \ln(\frac{z_{2,L}^*}{z_t})$ and $\tilde{\alpha} := \frac{2 |\alpha|}{\sigma^2}$. By some basic algebra, we see that

\[
  P(Z \mbox{ hits } z_{3,L}^* \mbox{ first} | Z(t)=z_t) = \frac{1 - (\frac{z_{2,L}^*}{z_t})^{\tilde{\alpha}} }{ (\frac{z_{3,L}^*}{z_t})^{\tilde{\alpha}} - (\frac{z_{2,L}^*}{z_t})^{\tilde{\alpha}}}.
\]

The split in the value function of the leader, see equation~\eqref{eq: leader_low_second_case}, appears because the leader does not know whether the demand process $Z_t$ will hit the innovation level of the follower, $z_{3,L}^*$, or the copying level of the follower, $z_{2,L}^*$, first. The splits have similar interpretations as those in Section~\ref{sec: F_L}.

%This expectation appears due to the leader not knowing whether the demand will reach the innovation level of the follower $z_{3,L}^*$ or the copying level $z_{2,L}^*$ first. The splits in the value function of the leader have similar interpretations as those in Section~\ref{sec: F_L}.

\medskip

The final case is where $(\frac{a_2}{a_1})^{\frac{\gamma}{\gamma - 1}} \geq \frac{K_2}{K_1}$. From Section~\ref{sec: F_L}, we know that in this case the follower-firm waits until the demand reaches the level $z_{3,L}^*$ and then innovates. By similar calculations as in the previous case, we find that

\begin{equation}
 \label{eq: leader_low_third_case}
 L_L(z_t) = \left\{
	\begin{array}{ll}
		\frac{\pi_l}{r-\alpha} z_t + \frac{\xi}{r-\alpha}z_t (1 - (\frac{z_t}{z_{3,L}^*})^{\gamma - 1} )   - \mc{I}  & \mbox{if } 0<z_t<z_{3,L}^* \\[\medskipamount]
		
		 \frac{ \pi_l }{r - \alpha}z_t - \mc{I} & \mbox{if } z_t \geq z_{3,L}^* \\[\medskipamount]
	\end{array}
\right.
\end{equation}

\subsubsection{The high profit case: $L_H(z_t)$}

Again, the leader's payoff function depends on the behavior of the follower. We have

\[
 \begin{array}{lll}
  L_H(z_t) &=& E^P[\int_t^{t+\tau^*} (\pi_h + \xi) Z_s e^{-r(s-t)} ds + \int_{\tau^* + t}^{\infty} \pi_h Z_s e^{-r(s-t)}ds - \mc{I}] \\[\smallskipamount]

   &=& E^P[\int_t^{\infty} \pi_h Z_s e^{-r(s-t)} ds] + E^P[\int_{t}^{\tau^* + t} \xi Z_s e^{-r(s-t)}ds] - \mc{I}. \\[\smallskipamount]
\end{array}
\]

The term $E^P[\int_{t}^{\tau^* + t} \xi Z_s e^{-r(s-t)}ds]$ equals zero if the follower stops immediately after the leader, and equals

\[
\frac{\xi}{r - \alpha}z_t (1 - (\frac{z_t}{z_H^*})^{(\gamma - 1)})
\]

for $z_t \in (0,z_H^*)$ (see Section \ref{lowprofcase}).

Thus, we have that

\begin{equation}
 \label{eq: leader_high}
 L_H(z_t) = \left\{
	\begin{array}{ll}
		\frac{\pi ^{h}}{r- \alpha}z_t + \frac{\xi}{r - \alpha}z_t(1 - (\frac{z_t}{z_H^*})^{(\gamma - 1)}) - \mc{I}  & \mbox{if } 0<z_t<z_H^* \\[\medskipamount]
		
		 \frac{ \pi_h}{r - \alpha}z_t - \mc{I} & \mbox{if } z_t \geq z_H^*. \\[\medskipamount]
	\end{array}
\right.
\end{equation}

\subsection{Combining the expressions}

To summarize, we have that the leader's ex ante expected payoff is

\begin{equation}
 \label{eq: final_leader_payoff}
 L(z_t) = \frac{1}{2}L_H(z_t) + \frac{1}{2}L_L(z_t),
\end{equation}
\noindent where we recall that $z^*_H < z^*_{1,L} < z^*_{3,L}$, $L_L(z_t)$ and $L_H(z_t)$ are as in equations~\eqref{eq: leader_low_third_case}-\eqref{eq: leader_high}.

Similarly, the follower's ex ante expected payoff is

\begin{equation}
 \label{eq: final_follower_payoff}
 F(z_t) = \frac{1}{2}F_H(z_t) + \frac{1}{2}F_L(z_t),
\end{equation}
\noindent where $F_L(z_t)$ and $F_H(z_t)$ are as in Section~\ref{sec: F_L} and equation~\eqref{eq: follower_high_2}.

By analyzing $L(z_t), F(z_t)$ and $C(z_t)$ (as functions of the demand $z_t$), we find that $C$ is a linear function satisfying $C(0) = -\mc{I}$. The leader payoff  function $L$ also satisfies $L(0)=-\mc{I}$. However, $L(z_t)$ is strictly concave on $(0, z_{1,L}^*)$, then equal to $C(z_t)$ for $z_t \in [z_{1,L}^*, z_{2,L}^*]$ and $z_t \geq z_{3,L}^*$. What happens to $L$ in the interval $[z_{2,L}^*, z_{3,L}^*]$ depends on the values of the parameters $\mc{I}, \theta, \alpha, r$. The follower's payoff $F$ satisfies $F(0)=-\mc{I}(1-\theta)$. It is strictly convex in the interval $(0, z_{1,L}^*)$. Also, it is greater than, but parallel to $C(z_t)$ (and therefore $L(z_t)$ as well) for $z_t \in [z_{1,L}^*, z_{2,L}^*]$. What happens in the interval $[z_{2,L}^*, z_{3,L}^*]$ depends on the parameters of the problem. For $z_t \geq z_{3,L}^*$, $F(z_t)$ is linear with a steeper slope than $C(z_t)$ and it is also larger than $C(z_t)$. An example of how the three payoff functions may look is shown in Figure \ref{fig:awesome_image} below.

\section{The stochastic game}
\label{Example_preemption}

In this section, we present a problem where the firms can choose when to invest, hence influencing whether they become the leader or the follower. Our aim is to illustrate that the timing of investment can become peculiar in a duopoly where there are benefits both of being leader (monopoly benefit) and follower (information benefit and reduced investment cost). The game we present is based on games discussed in Huisman et al.~\cite{Thijssen} and Chevalier-Roignant and Trigiorgis~\cite{TrigiorgisCompetitiveStrategies}.

\subsection{Defining the game}
\label{definingthegame}

We consider a duopoly, with the two firms denoted firm 1 and firm 2. For a subgame starting at time $t_0 \in [0, \infty)$, a simple strategy for firm $i$ ($i = 1,2$) is defined to be a function $ \alpha_i^{t_0} : [t_0, \infty) \times \Omega \rightarrow [0,1]$ such that
\begin{itemize}
\item $\forall \; \omega \in \Omega, \alpha_i^{t_0} (\cdot, \omega)$ is RCLL and right differentiable, and $\forall \; t \geq t_0, \alpha_i^{t_0}(t, \cdot)$ is $\mc{F}_t$-measurable.
\item If $\alpha_i^{t_0}(t,\omega) = 0$ for $t = \inf\{u \geq t_0 | \alpha_i^{t_0}(u;\omega) > 0\}$, then the right derivative of $\alpha_i^{t_0}(t ; \omega)$ is positive.
\end{itemize}
Furthermore, let  \[
\
\tau_i(t_0; \omega) := \left\{
	\begin{array}{ll}
		\infty  & \mbox{if } \alpha_i^{t_0}(t;\omega) = 0 \; \forall \mbox{ } t \geq t_0 \\[\medskipamount]
		
		 \inf\{t \geq t_0 | \alpha_i^{t_0}(t; \omega) > 0\} & \mbox{otherwise}, \\[\medskipamount]
	\end{array}
\right.
\]
and $\tau = \min\{\tau_1, \tau_2\}$.

The intuition of $\alpha_i^{t_0}$ can be interpreted in two steps. First, firm $i$ choosing $\alpha_i^{t_0}(t, \omega) = 0$ means that for the given scenario $\omega \in \Omega$, firm $i$ is not at all interested in becoming the leader firm at time $t$. If this holds for both firms in a neighborhood of $t$, then no investment will occur at this time. Second, firm $i$ choosing $\alpha_i^{t_0}(t, \omega) > 0$ means that firm $i$ definitely wants investment to occur at time $t$ for the given scenario, and is willing to be the leader in order for this to happen. However, the higher the level of $\alpha_i^{t_0}(t, \omega)$, the higher is the probability that firm $i$ actually becomes the leader given that also $\alpha_j^{t_0}(t, \omega) > 0$.

A closed loop strategy for firm $i$ is a collection of simple strategies 
\[
\{(\alpha_i^{t}(\cdot;\omega))_{0 \leq t < \infty}\} 
\] with the property that
\begin{itemize}
\item for all $u, v$ such that $ t \leq u \leq v < \infty,$ we have $\alpha_i^{t}(v; \omega) = \alpha_i^{u}(v; \omega)$ whenever $v= \inf\{s > t | Z_s = Z_v \}$.
\end{itemize}

The point of the above property is to ensure that the strategies of two different subgames coincide at all times $v$ that is a shared first hitting time of a new state $Z_v$ for both games. This ensures that strategies are consistent for the different subgames.

The expected discounted value of firm $i$ of the subgame starting at time $t_0$ given the simple strategies $s^{t_0}(\omega) = \{(\alpha_i^{t_0}(\cdot;\omega))_{i = 1,2}\}$ is
\begin{equation}
  \label{eq: V}
V_i(t_0,s^{t_0}(\omega)) = E_{t_0}[e^{-r\tau(t_0, \omega)}W_i(\tau(t_0, \omega), s^{t_0}(\omega))].
\end{equation}
Here, $E_{t_0}[\cdot] := E[\cdot|\mathcal{F}_{t_{0}}]$ and
 \[
W_i(\tau ,s^{t_0}(\omega)) = \left\{\begin{array}{ll}
L(Z_{\tau}(\omega)) & \text{if} \; \tau_j(t_0; \; \omega) > \tau_i(t_0; \; \omega),\\
F(Z_{\tau}(\omega)) & \text{if} \; \tau_i(t_0; \; \omega) > \tau_j(t_0; \; \omega).
\end{array} \right.
\]
Finally, if $\tau_i(t_0; \; \omega) = \tau_j(t_0; \; \omega)$, then
\[
W_i(\tau ,s^{t_0}(\omega)) = C(Z_{\tau}(\omega)) \mbox{ if } \; \alpha_i^{t_0}(\tau; \; \omega)) = \alpha_j^{t_0}(\tau; \; \omega)) = 1,
\]
\noindent while
\begin{equation}
\label{alphabetween}
\begin{array}{l}
W_i(t_0,s^{t_0}(\omega)) = \\
\left[\alpha_i^{t_0}(\tau; \; \omega)(1 - \alpha_j^{t_0}(\tau; \; \omega))L(Z_{\tau}(\omega)) +  \alpha_j^{t_0}(\tau; \; \omega)(1 - \alpha_i^{t_0}(\tau; \; \omega))F(Z_{\tau}(\omega)) \right.\\[\smallskipamount]
\left. + \alpha_i^{t_0}(\tau; \; \omega)\alpha_j^{t_0}(\tau; \; \omega)C(Z_{\tau}(\omega))\right] / [\alpha_i^{t_0}(\tau; \; \omega) + \alpha_j^{t_0}(\tau; \; \omega) - \alpha_i^{t_0}(\tau; \; \omega)\alpha_j^{t_0}(\tau; \; \omega)] \\[\smallskipamount]
\end{array}
\end{equation}
\noindent if  $2 > \alpha_i^{t_0}(\tau; \; \omega) + \alpha_j^{t_0}(\tau; \; \omega) > 0$ and
\begin{equation}
\label{alphazero}
W_i(t_0,s^{t_0}(\omega))= \frac{(\alpha_i^{t_0})'(\tau ; \; \omega)L(Z_{\tau}(\omega)) + (\alpha_j^{t_0})'(\tau ; \; \omega)F(Z_{\tau}(\omega))}{(\alpha_i^{t_0})'(\tau ; \; \omega) + (\alpha_j^{t_0})'(\tau ; \; \omega)} 
 \end{equation}
\noindent if $\alpha_i^{t_0}(\tau; \; \omega) = \alpha_j^{t_0}(\tau; \; \omega) = 0$. The functions $L, F$ and $C$ are as in Section \ref{sec: symmetric}.

For a given $\omega \in \Omega$ a tuple of simple strategies
\[
(s_i^{*})_{i=1,2}=(\alpha_i^{t_0}(t; \omega))_{i=1,2}
\]
\noindent is a Nash equilibrium for the subgame starting at time $t_0$ if for $i=1,2$ 
\[V_i(t_0, s_i^{*}, s_{-i}^{*}) \geq V_i(t_0, s_i, s_{-i}^{*})
\]
for all simple strategies $s_i$.

A pair of closed loop strategies $\{(\alpha_i^{t}(\cdot;\omega))_{0 \leq t < \infty}\}_{i= 1,2}$ is a subgame perfect equilibrium if for all $t \in [0, \infty)$ $(\alpha_i^{t}(\cdot; \omega))_{i = 1,2}$ is a Nash equilibrium.

\begin{figure}[!ht]
    \centering
    \includegraphics[width=0.5\textwidth]{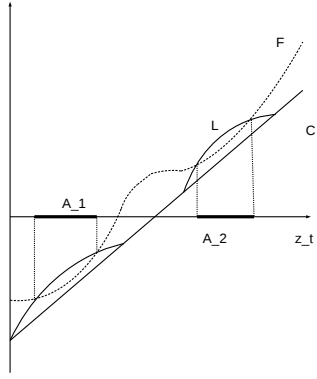}
    \caption{Illustration of the value functions for the leader $L$, the follower $F$ and the Cournot case $C$ (viewed as functions of the level of demand $z_t$), as well as the sets $A_1$ and $A_2$.}
    \label{fig:awesome_image}
\end{figure}

\subsection{The case with a high monopoly benefit}
\label{specialscenario}

We consider the case where $1< (\frac{a_2}{a_1})^{\frac{\gamma}{\gamma - 1}} <\frac{K_2}{K_1}$, so the function $F_L(z_t)$ is as in equation~\eqref{eq: follower_low_second_case}. Furthermore, we let the monopoly benefit $\xi$ be sufficiently large so that the following conditions hold
\begin{itemize}
\item there exists $z < z_{1,L}^{*}$ such that $L(z) > F(z)$,
\item  there exists $z \in (z_{2,L}^*, z_{3,L}^* )$ such that $L(z) > F(z)$.
\end{itemize}
This situation will always happen for $\xi$ sufficiently large, since we see from equation (\ref{eq: leader_low_second_case}) that $L(z)$ becomes arbitrarily large for $z < z_{1, L}^{*}$ and $z \in (z_{2,L}^*, z_{3,L}^* )$ as $\xi$ increases.

Thus, we are in the situation depicted in Figure \ref{fig:awesome_image}, where the follower and the leader value functions, $F(z)$ and $L(z)$, intersect twice. Let $A_1$ and $A_2$ denote the two intervals where the leader's payoff $L(z)$ exceeds the follower's $F(z)$. Also, recall that $z_{1,L}^*, z_{2,L}^*, z_{3,L}^*$ and $z_{H}^*$ are the investment thresholds of the follower in the low and high profit cases respectively, see Sections \ref{sec: F_L} and \ref{sec: F_H}.

Let $z_{\mbox{min}}^{A_1}$ be the infimum of $A_1$, and let similar definitions hold for $z_{\mbox{max}}^{A_1}, z_{\mbox{min}}^{A_2}$ and  $z_{\mbox{max}}^{A_2}$. Note that $z_{\mbox{min}}^{A_1}$ is the smallest level of demand such that $L$ is greater than or equal to $F$. We let $T_{z_{\min}^{A_1}}^{t}$ denote the stopping time for when the process is greater than or equal $z_{\mbox{min}}^{A_1}$ given we are at time $t$, with similar definitions for $z_{\mbox{max}}^{A_1}, z_{\mbox{min}}^{A_2}$ and $z_{\mbox{max}}^{A_2}$.

Under the assumptions above, we can prove the following result.

\begin{theorem}
\label{thm: theorema}
There is no pure symmetric subgame perfect Nash equilibrium in this setting. 

However, we have the following asymmetric equilibrium. For all $t \geq 0$, let
\[
\alpha_i^{t}(u) = \left\{
\begin{array}{ll}
\chi_{B_1}(Z_u)  & \text{if} \; u < T_{z_{\min}^{A_1}}^{t} \\
\frac{L(Z_u) - F(Z_u)}{L(Z_u) - C(Z_u)} & \text{if} \; T_{z_{\min}^{A_1}}^{t} \leq u < T_{z_{\max}^{A_1}}^{t} \\
\chi_{B_2}(Z_u) & \text{if} \; T_{z_{\max}^{A_1}}^{t} \leq u < T_{z_{\min}^{A_2}}^{t}\\
\frac{L(Z_u) - F(Z_u)}{L(Z_u) - C(Z_u)} & \text{if} \; T_{z_{\min}^{A_2}}^{t} \leq u < T_{z_{\max}^{A_2}}^{t} \\
\chi_{B_3}(Z_u) & \text{if} \; T_{z_{\max}^{A_2}}^{t} \leq u
\end{array} \right.
\]
and
\[
\alpha_j^{t}(u) = \left\{
\begin{array}{ll}
0 & \text{if} \; u < T_{z_{\min}^{A_1}}^{t} \\
\frac{L(Z_u) - F(Z_u)}{L(Z_u) - C(Z_u)} & \text{if} \; T_{z_{\min}^{A_1}}^{t} \leq u < T_{z_{\max}^{A_1}}^{t} \\
0 & \text{if} \; T_{z_{\max}^{A_1}}^{t} \leq u < T_{z_{\min}^{A_2}}^{t}\\
\frac{L(Z_u) - F(Z_u)}{L(Z_u) - C(Z_u)} & \text{if} \; T_{z_{\min}^{A_2}}^{t} \leq u < T_{z_{\max}^{A_2}}^{t} \\
0 & \text{if} \; T_{z_{\max}^{A_2}}^{t} \leq u,
\end{array} \right.
\]
where $i \neq j$. In the expression above, $\chi$ is the indicator function, and
\[
\begin{array}{l}
B_1 = \{y \in (0,z_{\min}^{A_1}) | L(y) \geq E[e^{-r T^{y}(z)}L(z)] \; \text{for all} \; z \in (0, z_{\min}^{A_1}] \}, \\
B_2 = \{y \in [z_{\max}^{A_1}, z_{\min}^{A_2}) | L(y) \geq E[e^{-r T^{y}(z)}L(z)] \; \text{for all} \; z \in [z_{\max}^{A_1}, z_{\min}^{A_2}] \}, \\
B_3 = \{y \in [z_{\max}^{A_2}, \infty) | L(y) \geq E[e^{-r T^{y}(z)}L(z)] \; \text{for all} \; z \in [z_{\max}^{A_2}, \infty) \},
\end{array}
\]
where $T^{y}(z)$ is the expected hitting time for the demand process  $Z_u$ to reach $z$, given that its current value is $y$.

\end{theorem}

\bigskip

\begin{proof}

The proof consists of two parts:

\begin{enumerate}
 \item
\emph{There is no subgame perfect symmetric Nash equilibrium:}

We start by showing that no symmetric equilibrium exist. For this, let $z_t$ be such that 
\begin{equation}
\label{Cineq}
C(z_t) > \sup\{L(z)| z \in (0, z_{\max}^{A_2}]\},
\end{equation}
and assume that there is a symmetric Nash equilibrium, i.e. simple strategies $\alpha_1, \alpha_2$ such that $\alpha_1^{t}(u) = \alpha_2^{t}(u) := \alpha^{t}(u)$ for all $u$.
\begin{itemize}
\item{If $\tau(t) = \infty$, i.e. $\alpha^{t}(u) = 0$ for all $u \geq t$, then both firms will regret their strategy. This follows since such a strategy $s^{t}(\omega)$ gives
\[
V_i(t,s^{t}(\omega)) = 0
\]
for $i = 1, 2$. If for instance firm 1 instead choose $\alpha_1^{t}(t) = 1$, then $\tau(t) = t$ and the payoff becomes $e^{-rt}L(z_t) > 0$ instead, which is better. Thus, we must have $\tau(t) < \infty$.}

\item{If $\tau(t) < \infty$, and if $F(Z_{\tau}) > L(Z_{\tau})$, then firm 1 will be better off by choosing a strategy such that $\tau_1 > \tau_2$, for instance by setting $\alpha_1^{t}(u) = 0$ for all $u$ such that $F(u) > L(u)$. This is because the former strategy gives payoff
\[
E_t[e^{-r\tau}W_1(\tau(t), s^{t}(\omega))],
\]
while any alternative strategy resulting in $\tau_1 > \tau_2$ gives payoff
\[
E_t[e^{-r\tau}F(Z_{\tau})]
\]
for firm 1. Since $F(Z_{\tau}) > L(Z_{\tau})$, it follows that $F(Z_{\tau}) > W_1(\tau(t), s^{t}(\omega))$, which makes such an alternative strategy strictly better.

If instead $F(Z_{\tau}) \leq L(Z_{\tau})$, then we must have $Z_{\tau} \leq z_{\max}^{A_2}$. However, it then follows that the payoff
\[
V_i(t,s^{t}(\omega)) < E_t[e^{-r \tau}C(z_t)] \leq e^{-r t}C(z_t),
\]
where we in the first inequality use the fact that (due to the description of $z_t$ from (\ref{Cineq}))
\[
C(z_t) > \sup\{L(z)| z \in (0, z_{\max}^{A_2}]\} \geq \sup\{F(z)| z \in (0, z_{\max}^{A_2}]\}.
\]
As a consequence firm 1 will be better off choosing $\alpha_1^{t}(t) = 1$, that is $\tau_1(t) = t$, as this gives the payoff $e^{-r t}C(z_t)$.}
\end{itemize}

Thus, we have shown that there does not exist any symmetric subgame perfect Nash equilibrium.

\item
\emph{The asymmetric strategy is an asymmetric subgame perfect Nash equilibrium:}

Next, we show that the described strategies is an asymmetric subgame perfect Nash equilibrium. Based on the form of the functions $\alpha_i^{t}, \alpha_j^{t}$, it is natural to classify the different time intervals into two categories:
\begin{itemize}
\item Category 1, where the firms' strategies are symmetric. That is, the time intervals $[T_{z_{\min}^{A_1}}^{t}, T_{z_{\max}^{A_1}}^{t}]$ and $[T_{z_{\min}^{A_2}}^{t}, T_{z_{\max}^{A_2}}^{t}]$.
\item Category 2, where the firms' strategies are not symmetric. That is, the time intervals $(0, T_{z_{\min}^{A_1}}^{t})$, $[T_{z_{\max}^{A_1}}^{t}, T_{z_{\min}^{A_2}}^{t})$ and $[T_{z_{\max}^{A_2}}^{t}, \infty)$.
\end{itemize}
 We now consider the two categories of time intervals separately.

\begin{itemize}
\item{Category 1. For $u \in [T_{z_{\min}^{A_1}}^{t}, T_{z_{\max}^{A_1}}^{t}] \cup [T_{z_{\min}^{A_2}}^{t}, T_{z_{\max}^{A_2}}^{t}]$, the strategies are chosen symmetrically in such a way that
\[
V_i(t, s^{t}(\omega)) = V_j(t, s^{t}(\omega)) = E_t[e^{-r u}F(Z_u))].
\]
Abbreviating and writing $s^{t}(\omega) = s = (s_i, s_j)$, we note that $s$ is chosen such that
\[
V_i(t, s_i', s_j) = V_i(t, s) = V_j(t, s) = V_j(t, s_i, s_j'),
\]
where $s_i', s_j'$ are any other strategy chosen by firm $i$ and firm $j$ respectively. This can be verified by evaluating the expression in (\ref{alphabetween}). That is, the strategy $s$ has the property that if only one firm changes strategy to any other strategy, the payoff of this firm will be unaltered. Thus, for any of the firms there is nothing to be attained by a change of strategy for time intervals belonging to Category 1.}

\item{Category 2. For $u \in (0, T_{z_{\min}^{A_1}}^{t}) \cup [T_{z_{\max}^{A_1}}^{t}, T_{z_{\min}^{A_2}}^{t}) \cup [T_{z_{\max}^{A_2}}^{t}, \infty)$, the strategy of firm $i$ is chosen in such a way that the firm will invest as long as investing and becoming leader is weakly better than waiting with the investment. Otherwise, firm $i$ will wait with the investment. Since firm $j$ is always waiting in time intervals belonging to Category 2, by the definition of the strategy, firm $i$ cannot be better off by any alternative strategy. Also, since the firms are ex ante identical, firm $j$ will not gain by choosing a strategy such that $\tau_j(u) < \tau_i(u)$ for any relevant $u$. And since we have $L < F$ in this interval, having $\alpha_j(u) = 0$ at times $u$ where $\alpha_i(u) > 0$ is an optimal strategy for firm $j$. We deduce that no firm gains by a change of strategy on time intervals belonging to Category 2.}
\end{itemize}

We conclude that the described strategy is a pure subgame perfect Nash equilibrium.
\end{enumerate}
\end{proof}

%We see that adding an extra investment choice which includes an information benefit to the follower greatly influences the outcome of the game, see for instance Fudenberg and Tirole~\cite{Fudenberg}.

For values $z_t \in A_1 \cup A_2$ the outcome is just as in the standard preemption game described in the literature, see e.g. Huisman et al.~\cite{Thijssen}. However, for all other starting values of $z_t$ the situation is very different. Due to the information spillover and lower price of investment, both firms have incentives to be the follower. Consequently, we get a situation where the firms need to weight the cost of investing before the other firm, versus the cost of waiting too long. One solution of this dilemma is for one firm to refuse to take any action at all, in which case the other firm behaves as in a monopoly setting, meaning the role as leader will be taken. How this will be done in practice, when both firms are equal, is an open question. In our described Nash equilibrium, firm $j$ has the characteristic of being more patient, while firm $i$ is more eager to invest quickly.

In our described equilibrium, the timing of the investment may become peculiar. If, for instance, $B_2 \neq [z_{\max}^{A_1}, z_{\min}^{A_2})$ (how $B_2$ looks depends on the parameters of the problem), then we have a scenario where
\begin{itemize}
\item no investment occurs for low values of $z_t$
\item investment occurs when $z_t \in [z_{\min}^{A_1}, z_{\max}^{A_1}]$
\item for some value $z_t \in (z_{\max}^{A_1}, z_{\min}^{A_2})$ no investment occurs
\item investment occurs for $z_t \in (z_{\min}^{A_2}, z_{\max}^{A_2}]$.
\end{itemize}
Thus a "vacuum" arises, where firms are happy to invest for lower and higher market demand levels, but not for intermediate market sizes. The explanation for such a vacuum is that in the set $A_1 \cup A_2$ the firms are preempted to invest immediately, since being leader is the most beneficial. Outside of this set the preemption effect vanishes, as firms are reluctant to invest since being the follower is more beneficial. In this region and in our described equilibrium, firm $i$ will behave as a monopolist having the payoff function $\min\{L(z_t), F(z_t)\}$ for all values $z_t$, given that firm $j$ is determined to become the follower.

\section{Conclusion}
\label{conclusion}
In this paper, we studied a duopoly investment model with uncertainty and two alternative investments. Being the first firm to invest gives a monopoly benefit for a period of time, while the follower firm attains information and cost reduction benefits. We have characterized the payoff functions for both the leader and follower firm. 

Furthermore, we have considered a game where the firms can choose when to invest, and hence influence whether they become the leader or the follower. The scenario we study leads to a more complex investment structure than what has been studied in previous literature. We show that no pure symmetric subgame perfect Nash equilibrium exists. However, an asymmetric equilibrium is characterized in Theorem~\ref{thm: theorema}. In this equilibrium, two disjoint intervals of market demand levels give rise to preemptive investment behavior of the firms, while the firms otherwise tend to be more reluctant to be the first investor.

\end{document}